\documentclass[twocolumn]{autart}    
 
\usepackage{graphicx}          
 
\usepackage{float}
\usepackage[colorlinks,linkcolor=blue]{hyperref}
\usepackage[latin1]{inputenc}
\usepackage{graphicx}
\usepackage{amssymb}
\usepackage{verbatim}
\usepackage{epsfig}
\usepackage[labelfont=bf]{caption}
\usepackage{enumerate}
\usepackage{subfig}
\usepackage{epstopdf}

\usepackage{stfloats}
\usepackage{amssymb}

\usepackage{graphicx}
\usepackage{bm}
\usepackage{color}
\usepackage{multirow}
\usepackage{mathrsfs}
\usepackage{lineno,hyperref}
\usepackage{lipsum}
\usepackage{stfloats}
\usepackage{titletoc}
\usepackage{appendix}
\usepackage{dsfont}
\usepackage{mathtools}

\def\begquo{\begin{quote}}
	\def\endquo{\end{quote}}
\def\begequarr{\begin{eqnarray}}
\def\endequarr{\end{eqnarray}}
\def\begequarrs{\begin{eqnarray*}}
	\def\endequarrs{\end{eqnarray*}}
\def\begarr{\begin{array}}
	\def\endarr{\end{array}}
\def\begequ{\begin{equation}}
\def\endequ{\end{equation}}
\def\lab{\label}
\def\begdes{\begin{description}}
	\def\enddes{\end{description}}
\def\begenu{\begin{enumerate}}
	\def\begite{\begin{itemize}}
		\def\endite{\end{itemize}}
	\def\endenu{\end{enumerate}}

\def\lef[{\left[\begin{array}}
	\def\rig]{\end{array}\right]}
\def\qed{\hfill$\Box \Box \Box$}
\def\begcen{\begin{center}}
	\def\endcen{\end{center}}
\def\begrem{\begin{remark}\rm}
	\def\endrem{\end{remark}}
\def\begdef{\begin{definition}}
	\def\enddef{\end{definition}}
	\def\begsubequ{\begin{subequations}}
\def\endsubequ{\end{subequations}}

\def\begmat#1{\begin{bmatrix}#1\end{bmatrix}}
\def\begali#1{\begin{align}{#1}\end{align}}
\def\begalis#1{\begin{align*}{#1}\end{align*}}

\def\caly{{\mathcal Y}}

\def\calm{{\mathcal M}}

\def\call{{\mathcal L}}

\def\bfp{{\bf p}}

\def\liminf{\lim_{t \to \infty}}

\def\L2e{{\cal L}_{2e}}

\def\rea{\mathbb{R}}

\def\adj{\mbox{adj}}
\def\col{\mbox{col}}

\def\rank{\mbox{rank}\;}

\def\phil{\phi^\mathtt{L}}

\usepackage{color}

\usepackage[prependcaption,colorinlistoftodos]{todonotes}

\begin{document}

\begin{frontmatter}

\title{Generalized Parameter Estimation-based Observers: Application to Power Systems and Chemical-Biological Reactors} 

\thanks[footnoteinfo]{This paper was not presented at any IFAC 
meeting. Corresponding author Nikolay~Nikolaev. Tel. +79213090016.}

\author[first,second]{Romeo~Ortega}
\ead{romeo.ortega@itam.mx},
\author[second]{Alexey~Bobtsov}
\ead{bobsov@mail.ru},
\author[second]{Nikolay~Nikolaev}
\ead{nikona@yandex.ru},
\author[third]{Johannes~Schiffer}
\ead{e-mail: schiffer@b-tu.de},
\author[fourth]{Denis~Dochain}
\ead{denis.dochain@uclouvain.be}

\address[first]{Departamento Acad\'{e}mico de Sistemas Digitales, ITAM, Ciudad de M\'exico, M\'{e}xico}
\address[second]{Department of Control
	Systems and Robotics, ITMO University, Kronverkskiy av. 49, Saint
	Petersburg, 197101, Russia}
\address[third]{Fachgebiet Regelungssysteme und Netzleitechnik, Brandenburgische Technische Universit\"{a}t Cottbus - Senftenberg}
\address[fourth]{CESAME, Universit\'{e} Catholique de Louvain (UCL), Avenue Georges Lema\^{i}tre , 4, B 1348, Louvain-la-Neuve, Belgium}


\begin{abstract}                          
In this paper we propose a new state observer design technique for nonlinear systems. It consists of an extension of the recently introduced parameter estimation-based observer, which is applicable for systems verifying a particular algebraic constraint. In contrast to the previous observer, the new one avoids the need of implementing an open loop integration that may stymie its practical application. We give two versions of this observer, one that ensures asymptotic convergence and the second one that achieves convergence in finite time. In both cases, the required excitation conditions are strictly weaker than the classical persistent of excitation assumption.   It is shown that the proposed technique is applicable to the practically important examples of multimachine power systems and chemical-biological reactors.
\end{abstract}

\end{frontmatter}
%
\section{Problem Formulation}
\label{sec1}
%
In this paper we are interested in the design of state observers for nonlinear control systems whose dynamics is described by
\begin{equation}
	\label{sys}
	\begin{aligned}
		\dot{x} & = f(x,u) \\
		y & = h(x,u),
	\end{aligned}
\end{equation}
where $x \in \rea^n$ is the systems state, $u \in \rea^m$ is the control signal and $y \in \rea^p$ are the {\em measurable} output signals. Similarly to all mappings in the paper, the mappings $f: \rea^n\times \rea^m  \to \rea^{n},\; h: \rea^{n} \times \rea^m \to \rea^{p}$ are assumed smooth. The problem is to design a dynamical system
{
\begin{equation}
	\begin{aligned}
		\dot{\chi} & = & F(\chi,y, u)\\
		\hat x & = & H(\chi,y,u)
		\label{dynobs}
	\end{aligned}
\end{equation}
with $\chi \in \rea^{n_\chi}$, such that {\em for all} initial conditions $x(0) \in \rea^n,\;\chi(0) \in \rea^{n_\chi}$, 
}
\begequ
\lab{obscon}
\liminf |\hat x(t)-x(t)|=0,
\endequ
where  $|\cdot|$ is the Euclidean norm. We are also interested in the case when the observer ensures {\em finite convergence time} (FCT), that is, when there exists $t_c \in [0,\infty)$ such that
\begequ
\lab{ftc}
\hat x(t)=x(t),\; \forall t \geq t_c.
\endequ
Following standard practice in observer theory \cite{BER} we assume that   $u$ is such that the state trajectories of \eqref{sys} are bounded.
Since the publication of the seminal paper  \cite{Luenberger1966TAC}, which dealt with linear time-invariant (LTI) systems, this problem has been extensively studied in the control literature. We refer the reader to  \cite{ASTetal,BER,BESbook,KHAbook} for a review of the literature. In this paper we propose an extension of the parameter estimation-based observer (PEBO) design technique reported in \cite{ORTetalscl}. The main novelty of PEBO is that it translates the task of state observation into an on-line  \emph{parameter estimation} problem.\\ 
The main features of the new observer design technique proposed in the paper, called generalized PEBO (GPEBO), are the following.\\
\begite
\item[(F1)] The key ``transformability into cascade form" condition of the original PEBO \cite[Assumption 1]{ORTetalscl} is {\em relaxed}, replacing it by a ``transformability into state-affine form" assumption discussed in \cite[Chapter 3]{BER}. 
\item[(F2)] We identify a class of systems for which the second key condition of PEBO \cite[Assumption 2]{ORTetalscl}---which relates with the, far from obvious, solution of the parameter estimation problem---is obviated. The class is identified   via a particular {\em algebraic constraint}.
\item[(F3)] It avoids the need of {\em open-loop integration} which stymies the practical application of this observer for systems subject to high noise environments---see  \cite[Remark R5]{ORTetalscl}.
\item[(F4)] Via the utilization of the fundamental matrix of an associated linear time-varying (LTV) system, the {\em signal excitation} needed to estimate the parameters is injected. 
\item[(F5)] Using the dynamic regressor extension and mixing (DREM) procedure \cite{ARAetaltac}, which is a novel, powerful, parameter estimation technique, we propose a variation of GPEBO achieving FCT, that is, for which \eqref{ftc} holds, under the weakest sufficient excitation assumption \cite{KRERIE}.\footnote{See \cite{GERetal} for an FCT version of DREM, \cite{ORTetalaut} for an interpretation as a Luenberger observer and \cite{ORTetaljpc,PYRetal} for two recent applications of DREM+PEBO techniques.} 
\item[(F6)]  It is proven that the conditions (F1) and (F2) are satisfied by the practically important case of multimachine power systems, while (F1) is verified by chemical-biological reactors. 
\endite
For the multimachine power systems we consider the well-know three-dimensional ``flux-decay" model of a large-scale power system \cite{KUN,SAUetal}, consisting of $N$ generators interconnected through a transmission network, which we assume to be lossy, that is, we explicitly take into account the presence of transfer conductances. We prove that, using the {\em measurements} of active and reactive power---which is a reasonable assumption given the current technology \cite{KUMPAL,SAUetal}---as well as the rotor angle at each generator, the application of GPEBO allows us to recover {\em the full state} of the system, even in the presence of lossy lines. To the best of the authors' knowledge, this is the first globally convergent solution to the problem. \\
For the reaction problem we consider the classical  dynamical model of the concentration components, {\em e.g.}, equation (1.43) in \cite[Section 1.5]{BASDOC}, which describes the behavior of a large class of chemical and bio-chemical reaction systems. We propose a state observer that, in contrast with the standard asymptotic observers \cite{BASDOC,DOCbook}, has a {\em tunable} convergence rate. Similarly to the case of power systems, using DREM, we can ensure FCT for the particular case when the reaction rates are linear in the unmeasurable states.\\
The remainder of the paper is organized as follows. In Section \ref{sec2}, to place in context the contributions of GPEBO, we briefly recall the basic principles of PEBO. In Section \ref{sec3} we give the main results. Section \ref{sec4} is devoted to some discussion. Section \ref{sec5} presents the application of the observer to an academic example and two practical problems. The paper is wrapped-up with concluding remarks in  Section \ref{sec6}. The proofs of the main propositions are given in appendices at the end of the paper.
%
\section{Review of PEBO and Introduction to GPEBO}
\lab{sec2}
%
To make the paper self-contained, in this section we briefly recall the underlying principle of our previous PEBO design \cite{ORTetalscl}. Then, with PEBO as the background, we highlight the main results of GPEBO, that extend its domain of applicability.
\subsection{Basic construction of PEBO}
\lab{subsec21}
As explained in the Introduction, the specificity of PEBO is that the problem of state observation is translated into a problem of {\em parameter estimation}, namely the initial conditions of the system \eqref{sys}. To achieve this objective, we consider in PEBO nonlinear systems of the form \eqref{sys} that can be transformed, via a change of coordinates, to a {\em cascade} form.  Let us assume for simplicity that the system is already given in this forms, namely that $\dot x = B(u,y)$ , where $B: \rea^m \times \rea^p \to \rea^{n}$ is some smooth mapping. In PEBO we do an open-loop integration of $B(u,y)$, that is we define
\begequ
\lab{dotxb}
\dot \xi = B(u,y),
\endequ
an operation that has well-known shortcomings---see  \cite[Remark R5]{ORTetalscl}---and make the observation that $\dot x =\dot \xi$, hence $x(t)=\xi(t)+\theta$ with $\theta:=x(0)=\xi(0)$. Then, construct the observed state as $\hat x =\xi+\hat \theta$, with $\hat \theta$ and estimate of the unknown vector $\theta$ using the information of $y$. Except for the case when $\theta$ enters linearly the task of generating a consistent estimate for $\theta$ is far from trivial. 
\subsection{New construction of GPEBO}
\lab{subsec22}
A first important difference of GPEBO is that we relax the assumption of transformability to a cascade form to transformability to an {\em affine-in-the-state} form
$$
\dot x = \Lambda(u,y) x + B(u,y)
$$
where $\Lambda: \rea^m \times \rea^p \to \rea^{n \times n}$. See \cite[Chapter 3]{BER} for a discussion on these normal forms and the existing observer designs for them.\\
In GPEBO the fragile step of open-loop integration of PEBO is replaced by the construction of a ``copy" of the system via another dynamical system 
$$
\dot \xi = \Lambda(u,y) \xi + B(u,y).
$$
avoiding the open-loop integration. The new idea introduced in GPEBO is to exploit the properties of the  {\em fundamental matrix} of an LTV system as follows. Let us define the signal 
\begequ
\lab{exxi0}
e=x-\xi, 
\endequ
whose dynamics is described by the LTV system
\begequ
\lab{dote0}
\dot e = A(t)e
\endequ
where $A(t):=\Lambda(u(t),y(t))$.  As shown in all textbooks of linear systems theory a property of LTV systems is that all solutions of \eqref{dote0} can be expressed as linear combinations of the columns of its fundamental matrix, which is the unique solution of the matrix equation
$$
\dot \Phi_A=A(t) \Phi_A,\;\Phi_A(0)=\Phi_A^0 \in \rea^{n \times n},
$$
with $\Phi_A^0$ full-rank, see \cite[Property 4.4]{RUGbook}. More precisely, 
$$
e(t)=\Phi_A(t) [\Phi_A^0]^{-1}e(0).
$$
Similarly to PEBO, in GPEBO we treat $e(0)$ as an unknown parameter $\theta:=e(0)$, that we try to {\em estimate}. Invoking \eqref{exxi0}, the observed state is then generated as
\begequ
\hat x = \xi + \Phi_A \hat \theta,
\endequ
where, to simplify notation and without loss of generality, we set $\Phi_A^0=I_n$, with $I_n$ the $n \times n$ identity matrix. The use of the fundamental matrix is the {\em key step} of GPEBO.\\
Another important advantage of GPEBO is that, if the output mapping $h(x,u)$ of \eqref{sys}, can also be expressed in an affine in the state form,\footnote{In our main result we consider a more general assumption, but here we use this simple one for the sake of clarity.} that is,
\begequ
\lab{hcd}
h(x,u)={\mathtt C}(u,y)x + {\mathtt D}(u,y),
\endequ
then it is possible to obtain a {\em linear regressor equation} (LRE) for the unknown vector $\theta$. Indeed, from the derivations above we get the LRE ${\mathtt y}=\psi \theta$
where we defined
\begalis{
	{\mathtt y} &:= y - {\mathtt D}(u,y)-{\mathtt C}(u,y)\xi\\
	\psi & :={\mathtt C}(u,y) \Phi_A.
}
This is also a fundamental feature since, as it is well-known \cite{LJUbook,SASBOD}, the design of parameter estimators for LRE is a well-understood problem. \\
A final advantage of GPEBO over PEBO pertains to the excitation conditions needed for parameter estimation. Notice that, if in PEBO the mapping $h(x,u)$ satisfies the assumption \eqref{hcd} we can also obtain a LRE ${\mathtt y}=\psi \theta$, but with $\psi={\mathtt C}(u,y)$. It is well-known that the convergence of all estimators is determined by the excitation of the regressor $\psi$. The presence of the additional term $\Phi_A$ in the regressor of GPEBO injects additional excitation. To appreciate this, consider the case when ${\mathtt C}(u,y)$ is {\em constant}. In that case it is impossible to estimate the parameter $\theta$ with the LRE of PEBO. 
%
\section{Main Results}
\lab{sec3}
%
The GPEBO designs are based on the following two propositions. For ease of presentation we consider the case where we are interested in observing {\em all} state variables. In many applications it is only necessary to reconstruct {\em some} of these state variables, a case that can be treated with slight modifications to these propositions. Also, we present first the version of GPEBO that ensures {\em asymptotic} convergence and then, in Proposition \ref{pro2}, the one ensuring FCT. The proofs of both propositions are given in Appendices \ref{appa} and \ref{appb}, respectively. 
\subsection{An asymptotically convergent GPEBO}
\label{subsec31}
%
\begin{prop}
	\label{pro1}
	\rm
	Consider the system \eqref{sys}. Assume there exist  mappings
	\begalis{
		& \phi: \rea^n \to \rea^{n},\;  \phil: \rea^{n} \times \rea^p \to \rea^{n},\;B: \rea^m \times \rea^p \to \rea^{n},\\
		&\Lambda: \rea^m \times \rea^p \to \rea^{n \times n},\;  L: \rea^{m}\times \rea^p \to \rea^{n \times n},\\
		&C: \rea^m \times \rea^p \to \rea^{n}
	}
	satisfying the following:
	\begite
	\item[(i)] The GPEBO partial differential equation (PDE)
	\begin{equation}
	\label{pde}
	\nabla \phi^\top(x) f(x,u) = \Lambda(u,h(x,u)) \phi(x) + B(u,h(x,u)),
	\end{equation}
	where $\nabla:=({\partial \over \partial x})^\top$.
	\item[(ii)] $\phil$ is a ``left inverse" of $\phi$, in the sense that it satisfies 
	\begequ
	\lab{phil}
	\phil(\phi(x),h(x,u))=x.
	\endequ
	\item[(iii)] The algebraic constraint
	\begequ
	\lab{algcon}
	L(u,h(x,u))\phi(x)=C(u,h(x,u)),
	\endequ
	is satisfied.
	\item[(iv)] For the given $u$, all solutions of the LTV system
	$$
	\dot z = \Lambda(u(t),y(t))z,
	$$
	with $y$ generated by \eqref{sys}, are {\em bounded}.
	\endite
	The GPEBO dynamics
	\begsubequ
	\lab{gpebodyn}
	\begali{
		\lab{dotxi}
		\dot{ \xi } & =  \Lambda(u,y)  \xi + B(u,y)\\
		\lab{dotphi}
		\dot \Phi_\Lambda &= \Lambda(u,y) \Phi_\Lambda,\;\Phi_\Lambda(0)=I_n\\
		\lab{doty}
		\dot Y &= - \lambda Y +  \lambda\Psi^\top [C(u,y) - L(u,y) \xi]\\
		\lab{dotome}
		\dot \Omega &=- \lambda \Omega +  \lambda\Phi_\Lambda\Phi_\Lambda^\top \\
		\dot {\hat \theta} &=-\gamma \Delta (\Delta \hat{\theta}-\caly),
		\lab{gpebo}
	}
	\endsubequ
	with $\lambda >0$ and $\gamma>0$, with the definitions
	\begsubequ
	\lab{gpebodyn1}
	\begali{
		\lab{Psi}
		\Psi&:=L(u,y) \Phi_\Lambda\\
		\lab{caly}
		\caly &:= \adj\{\Omega\}Y\\
		\Delta&:=\det\{\Omega\},
		\lab{del}
	}
	\endsubequ
	the state estimate
	\begequ
	\lab{hatx}
	\hat{x}  = \phil(\xi + \Phi\hat \theta,y),
	\endequ 
	ensures \eqref{obscon} with all signals bounded provided
	\begequ
	\lab{delcon}
	\Delta \notin \call_2.
	\endequ
	\qed
\end{prop}
\subsection{An GPEBO with FCT}
\label{subsec32}
%
A variation of GPEBO that ensures FCT is given in Proposition \ref{pro2}. To streamline its presentation we need the following sufficient excitation condition  \cite{KRERIE}.\footnote{This condition may be defined taking an initial time $t_0 > 0$ and integrating to $t_0+t_c$. Since we have fixed the initial time everywhere at zero we believe it is more appropriate to leave it like that.} 
\begin{assum}
	\lab{ass1}\em
	Fix a constant $\mu \in (0,1)$. There exists a time $t_c>0$ such that
	\begequ
	\lab{conint}
	\int_0^{t_c} \Delta^2(\tau) d\tau \geq - {1 \over \gamma} \ln(1-\mu).
	\endequ 
	\vspace{0.01cm}
\end{assum}
\begin{prop}\em
	\lab{pro2}
	Consider the system (\ref{sys}), verifying the conditions (i)-(iii) of Proposition \ref{pro1}. Fix $\gamma>0$ and $\mu \in (0,1)$. The state observer defined by \eqref{dotxi}-\eqref{gpebo}  and the state estimate
	\begali{
		\hat x &=  \phil\Big(\xi + \Phi_\Lambda {1 \over 1 - w_c}[\hat \theta - w_c \hat \theta(0)],y\Big),
		\lab{ftcgpebo}
	}
	with
	\begali{
		\dot w  &= -\gamma \Delta^2 w, \; w(0)=1,
		\lab{wwc}
	}
	and $w_c$ defined via the clipping function
	$$
	w_c = \left\{ \begin{array}{lcl} w & \;\mbox{if}\; & w < 1-\mu \\ 1-\mu & \;\mbox{if}\; & w \geq 1-\mu, \end{array} \right.,
	$$
	ensures \eqref{ftc} with all signals bounded provided  $\Delta$ verifies Assumption \ref{ass1}.
	\qed
\end{prop}
\section{Discussion}
\label{sec4}
%
\begenu[{\bf D1}]
\item The GPEBO PDE \eqref{pde} is a generalization of the PDEs that {are} imposed in the Kazantzis-Kravaris-Luenberger observer (KKLO), first presented in \cite{Kazantzis1998SCL} as an extension to nonlinear systems of Luenberger's observer, and further developed in  \cite{Andrieu2006SIAM}. Indeed, in KKLO the mapping $\Lambda(u,y)$ is a {\em constant}, Hurwitz matrix---see \cite{BERAND} for a recent extension to the non-autonomous case where the mapping $\phi$ depends on time (or the systems input). It also generalizes the PDE required in  PEBO where $\Lambda$ is {\em equal to zero}. 
\item As discussed in  \cite{ORTetalscl} and Section \ref{sec2}, a drawback of the original PEBO is that it involves an {\em open-loop} integration, namely
$$
\dot{ \xi } =  B(u,y),
$$
which stymies the practical application of PEBO in the presence of noise---see  \cite[Remark R5]{ORTetalscl}. Due to the presence of $\Lambda$ in the dynamics of $\xi$ given in \eqref{dotxi}, this difficulty is conspicuous by its absence in GPEBO. It should be pointed out that, using an alternative technique that relies on the Swapping Lemma \cite[Lemma 3.6.5]{SASBOD}, this shortcoming of PEBO has been overcome in \cite{PYRetalijc} for a class of electromechanical systems. 

\begin{figure}[h]
	\centering
	\includegraphics[width=0.45\linewidth]{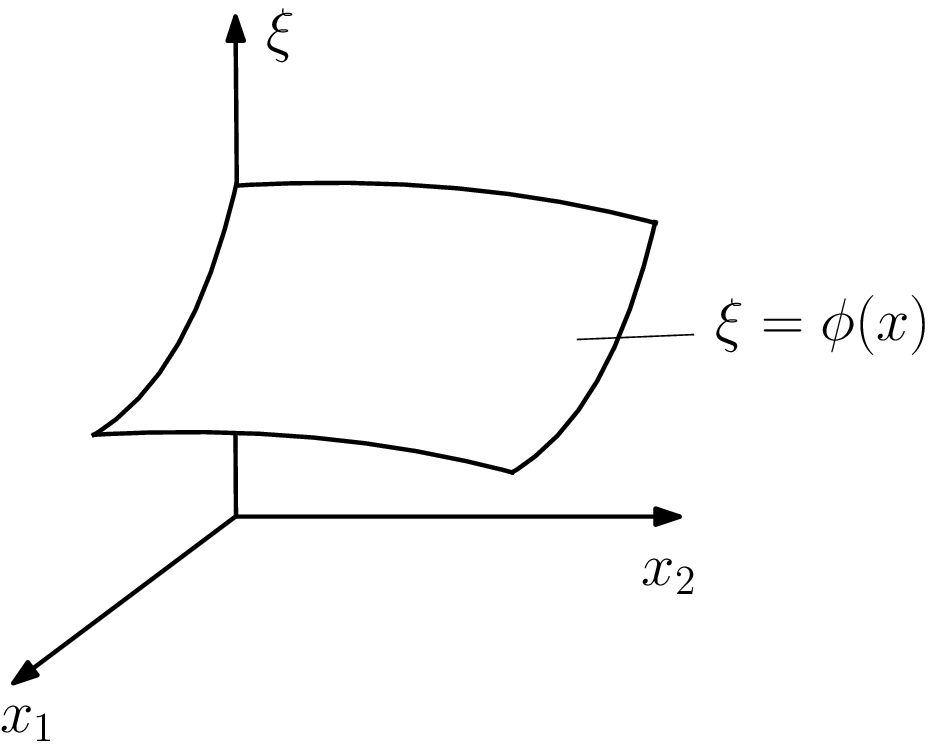}
	\includegraphics[width=0.45\linewidth]{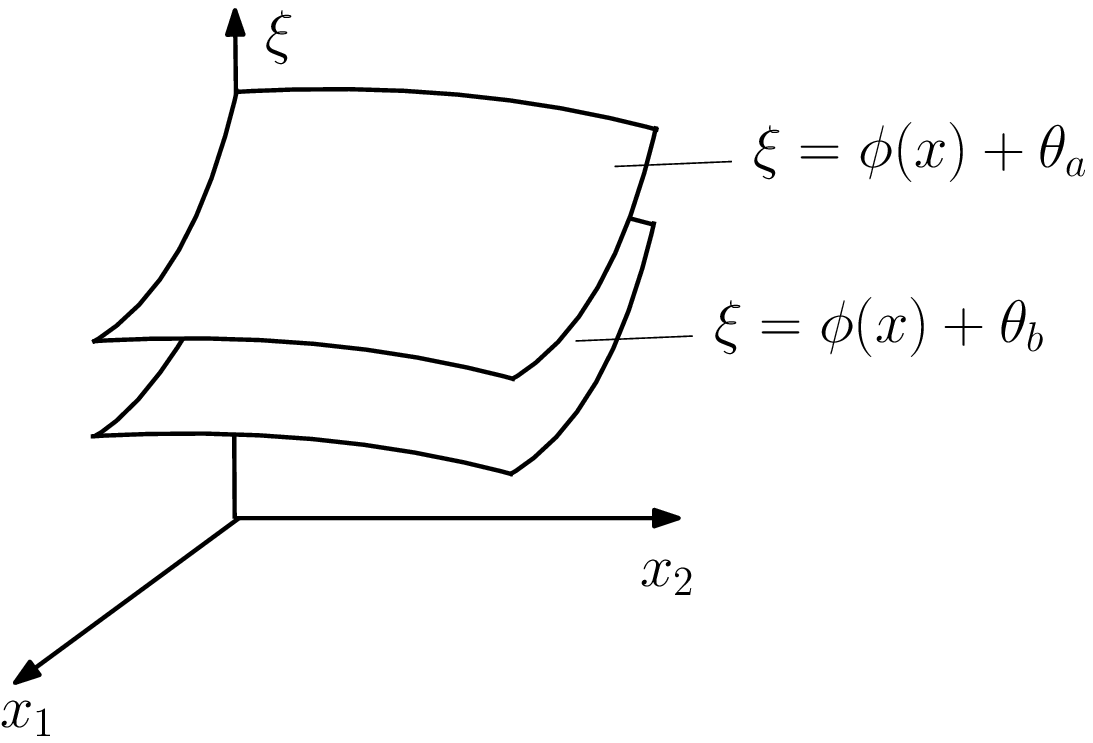}
	\caption{Geometric interpretation of KKLO and PEBO}
	\label{fig1}
\end{figure}

\item It is interesting to compare the KKLO with PEBO from the {\em geometric} viewpoint. The former generates an {\em attractive and invariant manifold}
$$
\calm:=\{ (x,\xi) \in \rea^n \times \rea^{n}\;|\; {\xi = \phi(x)}\},
$$
and the state is reconstructed, via $\phil$, with  $\xi$. On the other hand, PEBO generates an {\em invariant foliation}
{ 
	$$
	\calm_\theta:=\{ (x,\xi) \in \rea^n \times \rea^{n}|{\xi = \phi(x)+\theta}, \theta \in \rea^{n}\},
	$$
}
that is, the sublevel sets of the function $F(\xi,x):=\xi - \phi(x)$. To reconstruct the state---again via $\phil$---it is necessary to identify the leaf of $\calm_\theta$ via the estimation of $\theta$. See Fig. \ref{fig1}. See also \cite{YIetal} where it is proposed to combine PEBO and KKLO to extend the realm of application of these observers.

\item Imposing the algebraic constraint (ii) of Proposition \ref{pro1} is, clearly, a strong assumption. It is interesting that---as shown in Section \ref{sec4}---it is satisfied for the, practically relevant, power systems example. See also  \cite{PYRetalijc} where similar constraints are shown to be satisfied by a class of electromechanical systems and  \cite{PYRetal} for a significant extension, to the case of {\em adaptive} state observers---that is, systems with uncertain parameters and unmeasurable states---is reported.  

\item The version of DREM utilized in Proposition \ref{pro1} uses the dynamic extension proposed by \cite{KRE}. As discussed in \cite{ORTetaltac} other versions of DREM, with different convergence properties, are also possible. We have opted for this variation for the sake of simplicity.

\item The conditions $\Delta \notin \call_2$ and Assumption \ref{ass1} are, evidently, excitation conditions necessary to ensure convergence of the parameter estimators. Clearly, this kind of assumptions are unavoidable in the problem of state (or parameter) estimation. It is interesting that, as shown in \cite{ORTetaltac},  these conditions are {\em strictly weaker} than the usual persistent of excitation assumption imposed in standard parameter estimation schemes \cite[Theorem 2.5.1]{SASBOD}. 

\item It is possible to obviate the parameter estimation step of PEBO designing a KKLO-like observer. Indeed, under assumptions (i)-(iii) of Proposition \ref{pro1} the observer of $\phi$
\begin{footnotesize}
	$$
	\dot{\hat \phi}=  \Lambda(u,y)  \hat \phi + B(u,y)+\gamma  L^\top (u,y)[C(u,y) - L(u,y) \hat \phi]
	$$
\end{footnotesize}
verifies the error model
$$
\dot{\tilde \phi}=  [\Lambda(u,y) - \gamma  L^\top (u,y)L(u,y)] \tilde \phi.
$$
where $\tilde \phi:=\hat \phi- \phi$. However, some additional assumptions have to be imposed to the mappings $\Lambda$ and $L$ to ensure asymptotic stability of this LTV system.
\endenu
\section{Applications}
\label{sec5}
In this section we illustrate with an academic example and two physical systems the applicability of the proposed GPEBO. Towards this end, we identify all the mappings required to verify some (or all) of the conditions of Proposition \ref{pro1}. 
\subsection{An academic example}
\label{subsec51}
In \cite{BERAND} the problem of state observation of the following system is considered
\begali{
	\nonumber
	\dot x_1&=x_2^3\\
	\nonumber
	\dot x_2&=-x_1\\
	\lab{pausys}
	y&=x_1.
}
\subsubsection{Solution via PEBO+DREM}
\label{subsubsec511}
The proposition below shows that this problem can be trivially---and robustly---solved using the classical PEBO+DREM approach. Indeed, with $B=-y$ we see that the unknown part of the state, {\em e.g.} $x_2$ can be written in the form \eqref{dotxb}. Hence, following the PEBO procedure, we define
\begequ
\lab{dotxipau}
\dot \xi=-y. 
\endequ

\begin{prop}\em
	\lab{pro3}
	Consider the system (\ref{pausys}). Define  the state estimate
	$$
	\hat{x}_2  = \xi + \hat \theta,
	$$ 
	with \eqref{dotxipau} and the {\em scalar} parameter estimator
	$$
	\dot {\hat \theta} =-\gamma \Delta (\Delta \hat{\theta}-\caly_1),\;\gamma>0,
	$$
	where we define
	\begalis{
		Y&={\lambda \bfp \over \bfp+\lambda}[y] -{\lambda  \over \bfp+\lambda}[\xi^3],\;\\
		\Omega &={\lambda  \over \bfp+\lambda}[\phi \phi^\top ], \; 
		\phi=\begmat{{3\lambda  \over \bfp+\lambda}[\xi^2] \\ {3\lambda  \over \bfp+\lambda}[\xi] \\ {\lambda  \over \bfp+\lambda}[u_{-1}(t)]},\; \\ \caly &= \adj\{\Omega\}{\lambda  \over \bfp+\lambda}[\phi Y],\;\Delta=\det\{\Omega\},
	}
	where $\bfp:={d \over dt}$ and $u_{-1}(t)$  is a step signal happening at $t=0$.\footnote{We use the step signal convention to indicate that $\phi_3$ is the solution of the differential equation $\dot{\phi_3} = -\lambda \phi_3 + \lambda$, with $\phi_3(0) = 0$.} Then, \eqref{obscon} holds provided \eqref{delcon} is verified. 
\end{prop}

\begin{pf}
	Clearly $x_2=\xi+\theta$, where $\theta:=x_2(0)-\xi(0)$. Replacing in \eqref{pausys}, and developing the cubic power yields
	$$
	\dot x_1=\xi^3+3\xi^2 \theta+3 \xi \theta^2+\theta^3.
	$$
	Applying the filter ${\lambda  \over \bfp+\lambda}$, and using the definitions of $Y$ and $\phi$ above yields 
	$$
	Y=\phi^\top  \Theta,\; \Theta:=\begmat{\theta \\ \theta^2 \\ \theta^3}.
	$$ 
	Multiplying the equation above by $\phi$, applying again the filter ${\lambda  \over \bfp+\lambda}$, and multiplying by $\adj\{\Omega\}$ yields $\caly = \Delta \Theta$. The proof is completed replacing $\caly_1 = \Delta \theta$ in  the parameter estimator to get the error equation \eqref{dottilthe}.   
\end{pf}
\subsubsection{Simulations}
\label{subsubsec512}
\begin{figure}[h]
	\centering
	\includegraphics[width=1\linewidth]{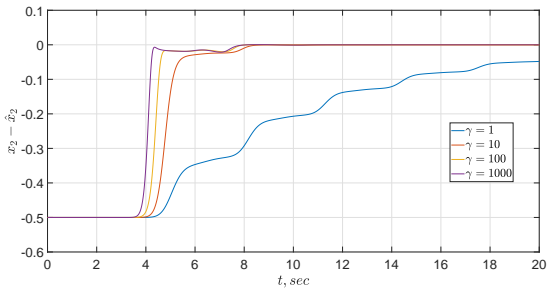}
	\caption{Transients of $x_2-\hat{x}_2$}
	\label{fig2}
\end{figure}

In Fig. \ref{fig2} we show the simulation results of the observer of Proposition \ref{pro3} with $x_1(0)=1$, $x_2(0)=0$, $\lambda=1$, $\hat{\theta}(0)=0.5$ and all filters initial conditions (ICs) zero. Notice that during the first $3$ seconds the estimates are ``frozen". This is due to the fact that, because of our choice of the observer ICs, the matrix $\Omega$ is rank deficient. Also, as expected, the rate of convergence is improved increasing the adaptation gain $\gamma$. These transients should be compared with the ones shown in  Fig. 1 of \cite{BERAND}, which are generated with a far more complicated KKLO that, moreover, needs to estimate also $x_1$. On one hand, the convergence time of the GPEBO design is significantly faster than the KKLO. On the other hand, our response is quite smooth while the one of the KKLO exhibits large oscillations.  
\subsection{Multimachine power systems}
\label{subsec52}
The dynamical model of the $i$--th generator of  $n$ interconnected machines can be described using the classical third order model\footnote{To simplify the notation, whenever clear from the context, the qualifier ``$i\in \bar n$" will be omitted in the sequel.}  \cite{KUN,SAUetal}
\begin{equation}
\begin{split}
\dot{\delta}_i&=\omega_i \\
M_i\dot{\omega}_i&=-D_{mi}\omega_i+\omega_0(P_{mi}-P_{ei})\\
\tau_i\dot{E}_i&=-E_i-(x_{di}-x'_{di})I_{di}+E_{fi}+\nu_i,\;\;\;\\
&i\in \bar n:=\{1,...,n\},
\label{syspow}
\end{split}
\end{equation}  
where the state variables  are the rotor angle {$\delta_i\in \mathbb R$ ,}
$\mbox{rad}$, the speed deviation $\omega_i\in\mathbb{R}$ in $\mbox{rad/sec}$ and the generator
quadrature internal voltage $E_i\in\mathbb{R}_+$, $I_{di}$ is the $d$ axis current, $P_{ei}$ is the electromagnetic power, the voltages $E_{fi}$ and
$\nu_i$ are the constant voltage component applied to the field winding, and the control voltage input, respectively. $D_{mi}$, $M_i$, $P_{mi}$,  $\tau_i$, $\omega_0$, $x_{di}$ and $x'_{di}$ are positive parameters.\\ 
The active power $P_{ei}$ and reactive power $Q_{ei}$ are defined as
\begin{equation}
P_{ei}=E_i I_{qi}, \;
Q_{ei}=E_i I_{di},
\label{peqe}
\end{equation}
where $I_{qi}$ is the $q$ axis current.\\
These currents establish the connections between the machines and are given by
\begin{equation} 
\begin{split}
I_{qi}&=G_{mii}E_{i}+ \sum_{j=1,j\neq i}^{n} E_{j} Y_{ij} \sin(\delta_{ij}+\alpha_{ij})\\
I_{di}&=-B_{mii}E_{i}-\sum_{j=1,j\neq i}^{n} E_{j} Y_{ij} \cos(\delta_{ij}+\alpha_{ij}),
\label{idiq}
\end{split}
\end{equation}
where we defined $\delta_{ij}:=\delta_i-\delta_j$ and the constants $Y_{ij}=Y_{ji}$ and $\alpha_{ij}=\alpha_{ji}$ are the admittance magnitude and admittance angle of the power line connecting nodes $i$ and $j$, respectively. Furthermore, $G_{mii}$ is the shunt conductance and $B_{mii}$ the shunt susceptance at node $i$. Finally, combining (\ref{syspow}), \eqref{peqe} and (\ref{idiq}) results in the well-known compact form 
\begin{equation} 
\begin{split}
\dot{\delta}_i &=  \omega_i\\
\dot{\omega}_i &=  -D_i\omega_i+P_i-d_i\Big[G_{mii}E_i^2 \\
&- E_i\sum_{j=1,j\neq i}^n E_j Y_{ij}\sin(\delta_{ij}+\alpha_{ij})\Big]\\
\dot{E}_i &=  -a_iE_i+b_i\sum_{j=1,j\neq i}^n E_j Y_{ij}\cos(\delta_{ij}+\alpha_{ij})+ u_i,
\label{syscom}
\end{split}
\end{equation}
where we have defined the signal
$$
u_i:= \frac{1}{\tau_i}(E_{f_i}+\nu_i)
$$
and the positive constants 
\begalis{
	D_i&:=\frac{D_{mi}}{M_i},\; P_i:=d_i P_{mi},\; d_i:=\frac{\omega_0}{M_i}\\
	a_i&:=\frac{1}{\tau_i}[1-(x_{di}-x'_{di})B_{mii}],\;b_i:=\frac{1}{\tau_i}(x_{di}-x'_{di}).
}
To formulate the observer problem we consider that all parameters are known, and make the following assumption on the available {\em measurements}.
\begin{assum}
	\label{ass2} \em
	The signals {$u_i$,} $\delta_i$, $P_{ei}$ and $Q_{ei}$ of all generating units are {\em measurable}.
\end{assum}
It is fair to say that the assumption of knowledge of $\delta_i$ is far from realistic.
\subsubsection{Verifying the conditions of Proposition \ref{pro1}}
\label{subsubsec521}
We make the following observation. Using \eqref{peqe} and \eqref{idiq}, the rotor speed dynamics \eqref{syscom} may be written as
$$
\dot{\omega}_i =  -D_i\omega_i+P_i-d_i P_{ei}.
$$
Considering that $P_{ei}$ is measurable, while $P_i$, $D_i$ and $d_i$ are known positive constants, the design of an observer for this system is trivial. For instance, 
\begin{equation}
	\begin{split}
			\dot\xi_{\omega_i} &=   -D_i\hat \omega_i+P_i-d_iP_{ei} - k_{\omega i} \hat{\omega}_i\\
		\hat \omega_i &=\xi_{\omega_i}+ k_{\omega i} \delta_i,\;k_{\omega i}>0,
	\end{split}
\label{hat_om}
\end{equation}
yields the LTI, asymptotically stable error dynamics 
$$
\dot {\tilde \omega}_i =  -(D_i+k_{\omega i})\tilde \omega_i.
$$
Therefore, we concentrate in the estimation of the voltages $E_i$. Its dynamics may be written as
\begequ
\lab{dote}
\dot E = \Lambda(\delta)E+u.
\endequ 
where $E:=\col(E_1,\dots,E_n)$, $\delta:=\col(\delta_1,\dots,\delta_n)$, and we defined matrix
\begali{
	\Lambda(\delta):=(\Lambda_1(\delta) \; \Lambda_2(\delta) \; \dots \; \Lambda_n(\delta)),
			\lab{a}
} 
where
\begali{
\nonumber
\Lambda_1(\delta)&:= \begmat{-a_1 \\ b_2 Y_{21}\cos (\delta_{21}+\alpha_{21}) \\ b_n Y_{n1}\cos (\delta_{n1}+\alpha_{n1})},\\
\nonumber
\Lambda_2(\delta)&:=\begmat{b_1 Y_{12}\cos (\delta_{12}+\alpha_{12}) \\ -a_2 \\ b_n Y_{n2}\cos (\delta_{n2}+\alpha_{n2})},\\
\nonumber
\Lambda_n(\delta)&:=\begmat{b_1 Y_{1n}\cos (\delta_{1n}+\alpha_{1n}) \\ b_2 Y_{2n}\cos (\delta_{2n}+\alpha_{2n}) \\ -a_n}
}
and we recall that $\delta$ is {\em measurable}. The remaining mappings of (i) and (ii) of Proposition \ref{pro1} are given as $\phi=E$ and $B=u$. The following simple lemma defines the mappings $L$ and $C$ that satisfy \eqref{algcon}.
\begin{lem}\em 
	\lab{lem1}
	There exists a {\em measurable} matrix $L(P_e,Q_e,\delta)\in \rea^{n \times n}$ such that
	\begequ
	\lab{ce}
	LE=0.
	\endequ
Consequently, selecting $C=0$, \eqref{algcon} is satisfied 
\end{lem} 
\begin{pf}
	From \eqref{peqe} we have that
	$$
	P_eI_d - Q_eI_q=0.
	$$
	Clearly, the equations \eqref{idiq}---which are linearly dependent on $E$---may be written in the compact form
	\begequ
	\lab{idiqcom}
	I_q=S(\delta)E,\;I_d=T(\delta)E,
	\endequ 
	for some suitably defined $n \times n$ matrices $S(\delta),\;T(\delta)$. The proof is completed by replacing \eqref{idiqcom} in the identity above and defining
	$$
	L(P_e,Q_e,\delta):=\begmat{P_{e1}T_1^\top (\delta)-Q_{e1}S_1^\top (\delta) \\ \vdots \\ P_{en}T_n^\top (\delta)-Q_{en}S_n^\top (\delta)},
	$$
	where $T_i^\top (\delta),\;S_i^\top (\delta)$ are the rows of the matrices $T(\delta)$ and $S(\delta)$, respectively.
\end{pf}
This lemma completes the verification of all the conditions of Proposition \ref{pro1}.
\subsubsection{Simulations}
\label{subsubsec522}
For simulation we use the two-machine system considered in \cite{ORTGALAST}. The dynamics of the system result in the sixth-order model
\begin{footnotesize}
	\begin{align}
	\label{2masch}
	\begin{cases}
	\dot{\delta}_1 &=\omega_1, \\
	\dot{\omega}_1 &=-D_1 \omega_1 + P_1 - G_{11} E_1^2 - Y_{12} E_1 E_2 \sin(\delta_{12}+\alpha_{12})\\
	\dot{E}_1 &= -a_1 E_1 + b_1 E_2 \cos(\delta_{12}+\alpha_{12}) + E_{f_1} + \nu_1;\\
	\dot{\delta}_2 &=\omega_2, \\
	\dot{\omega}_2 &=-D_2 \omega_2 + P_2 - G_{22} E_2^2 + Y_{21} E_1 E_2 \sin(\delta_{12}+\alpha_{12})\\
	\dot{E}_2 &= -a_2 E_2 + b_2 E_1 \cos(\delta_{21}+\alpha_{21}) + E_{f_2} + \nu_2,
	\end{cases}
	\end{align}
\end{footnotesize}
with the current equations defined as 
\begin{align}
I_{q1}&=G_{11}E_{1}+ E_{2} Y_{12} \sin(\delta_{12}+\alpha_{12})\\
I_{d1}&=-B_{11}E_{1}-E_{2} Y_{12} \cos(\delta_{12}+\alpha_{12})\\
I_{q2}&=G_{22}E_{2}+ E_{1} Y_{21} \sin(\delta_{21}+\alpha_{21})\\
I_{d2}&=-B_{22}E_{2}-E_{1} Y_{21} \cos(\delta_{21}+\alpha_{21}).
\end{align} 
In this case we have that
\begalis{
	A(t)&=\begmat{ -a_1 & b_1  \cos(\delta_{12}(t)+\alpha_{12}) \\ b_2  \cos(\delta_{21}(t)+\alpha_{21})& -a_2 }\\
	S(\delta)&= \begmat{ G_{11} & Y_{12} \sin(\delta_{12}+\alpha_{12})\\  Y_{21} \sin(\delta_{21}+\alpha_{21}) & G_{22}}\\
	T(\delta)&= \begmat{- B_{11} & - Y_{12} \cos(\delta_{12}+\alpha_{12})\\  - Y_{21} \cos(\delta_{21}+\alpha_{21}) & - B_{22}}.
}
For the observer design we selected the simplest filter
$$
F(p)=\begmat{1 & 0 \\ \frac{k}{p+k} & 0},
$$
with $p:=\frac{d}{dt}$ and $k>0$. The parameters of the model \eqref{2masch} are taken from \cite{ORTGALAST} and are given in the table given in Appendix \ref{appc}.\\
{Simulation results are presented in Fig. \ref{fig:pic2}-Fig. \ref{fig:pic6}. {Fig. \ref{fig:pic2} and Fig. \ref{fig:pic4} show the observation errors for  the open loop observer (OLO) \eqref{dotxi}, and for DREM for different adaptation gains and for FCT-DREM. For simulation we used $\lambda=1$ in \eqref{doty} and \eqref{dotome}. Simulation results for FCT-DREM are preseted for $\gamma=10^7$ in \eqref{gpebo} and $\mu=0.1$ for computation $\omega_c$ in \eqref{ftcgpebo}.} To test the robustness of the desing a $30\%$ load change was introduced at $t=10$ sec, whose effect is impercebtible. {Fig. \ref{fig:pic5} and Fig. \ref{fig:pic6} show the observation errors for rotor speed observer \eqref{hat_om} for first and second generator for different values of $k_{\omega i}$ in \eqref{hat_om}.} 
\begin{figure}[hbtp]
	\begin{center}
		\includegraphics[width=1 \linewidth]{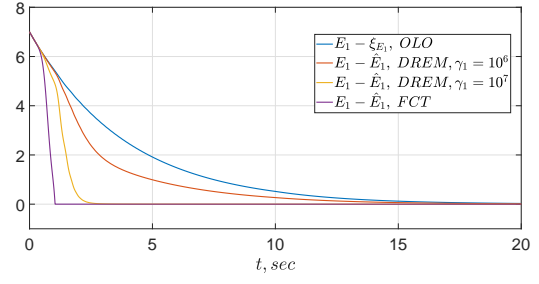}
				\caption{Transients of the first voltage observation error for the OLO \eqref{dotxi}, DREM and FCT-DREM observers with a 30\% load change a $t=10$ sec
		}
		\label{fig:pic2}
	\end{center}
\end{figure}
\begin{figure}[hbtp]
	\begin{center}
	\includegraphics[width=1 \linewidth]{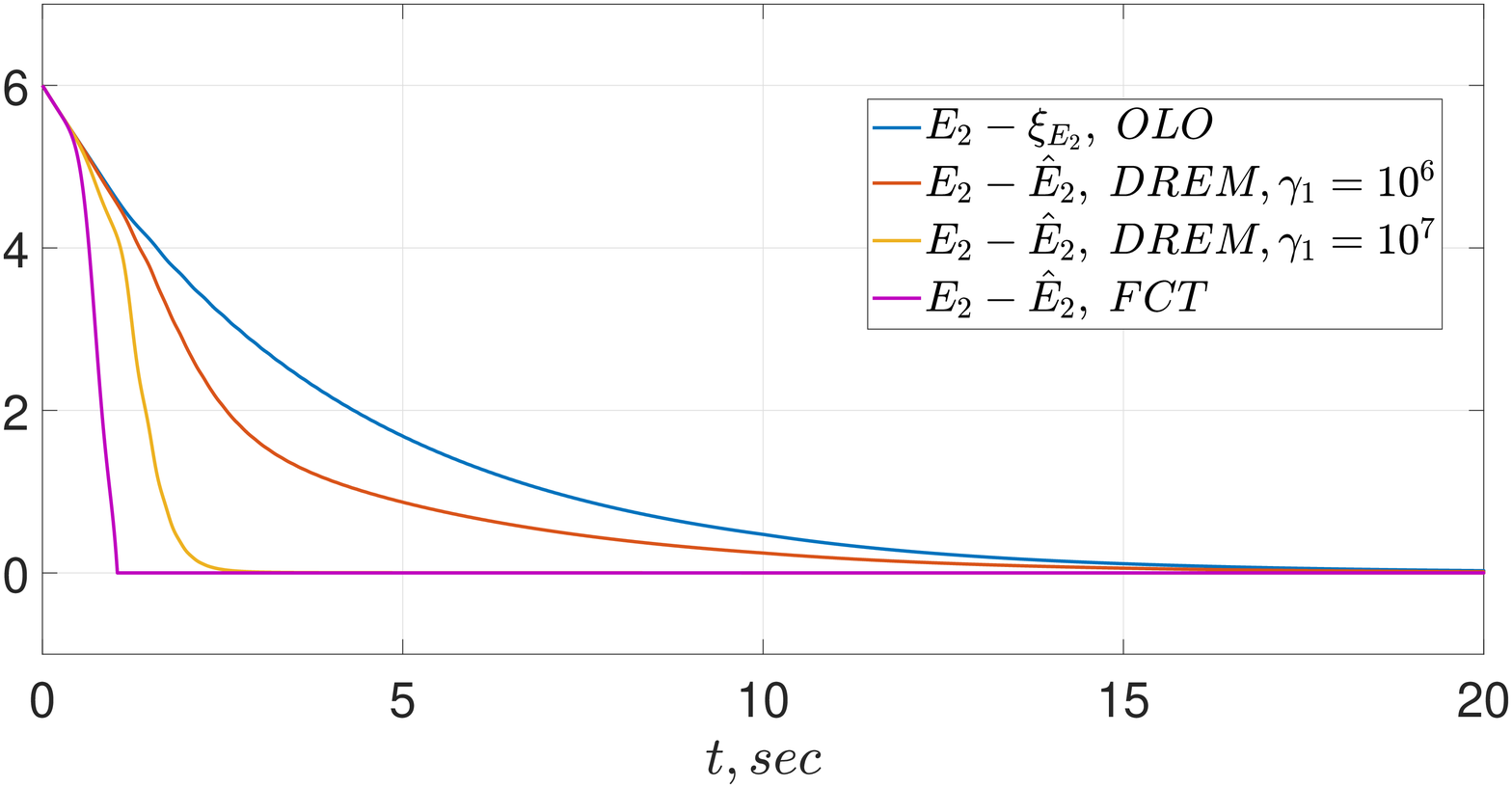}    
		\caption{Transients of the second voltage observation error for the OLO \eqref{dotxi}, DREM and FCT-DREM observers with a 30\% load change at $t=10$ sec
		}
		\label{fig:pic4}
	\end{center}
\end{figure}
\begin{figure}[hbtp]
	\begin{center}
		\includegraphics[width=1 \linewidth]{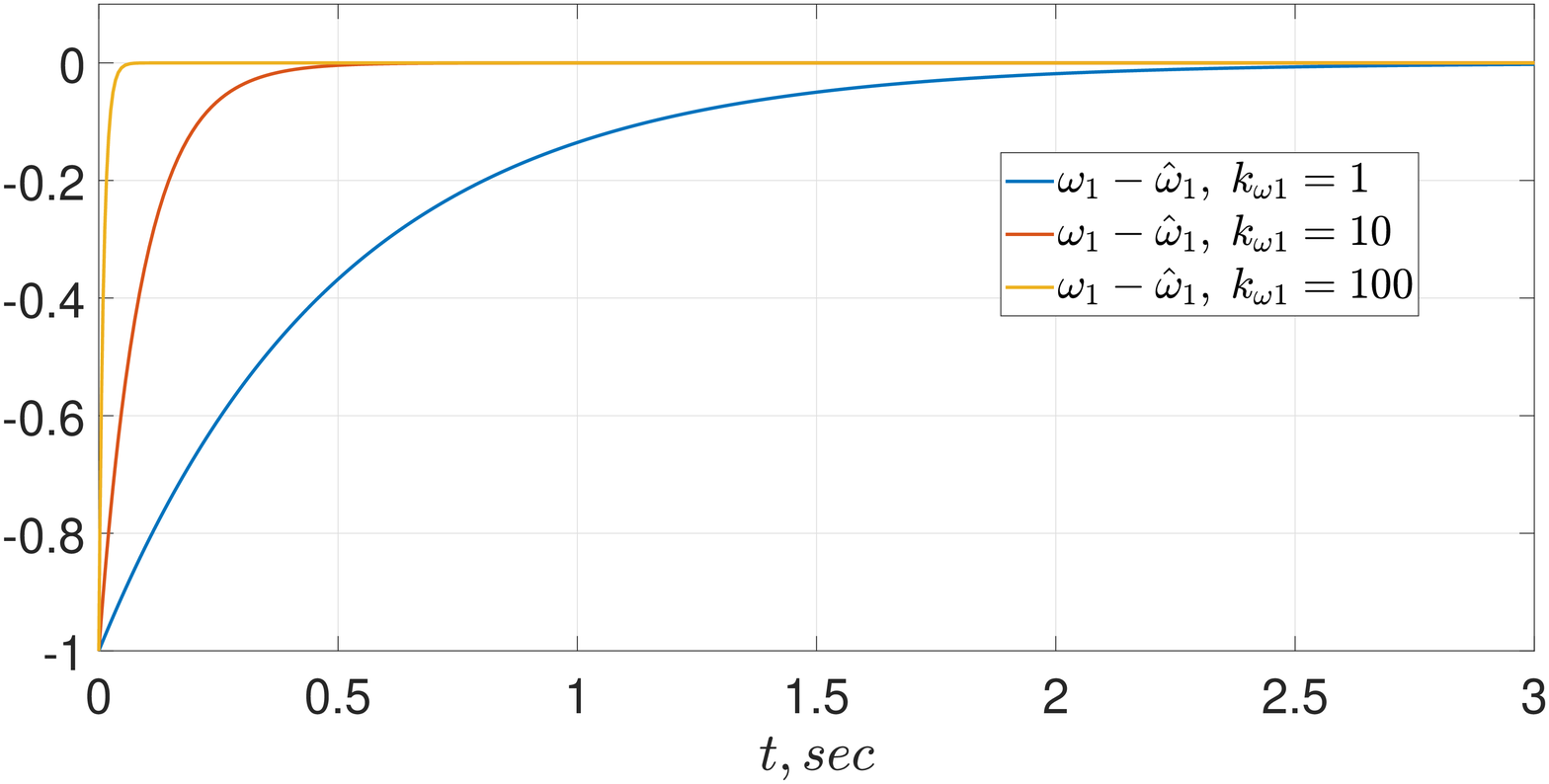}    
		\caption{Transients of the first speed observation error for the observer \eqref{hat_om} for different values of $k_{w1}$
		}
		\label{fig:pic5}
	\end{center}
\end{figure}
\begin{figure}[hbtp]
	\begin{center}
		\includegraphics[width=1 \linewidth]{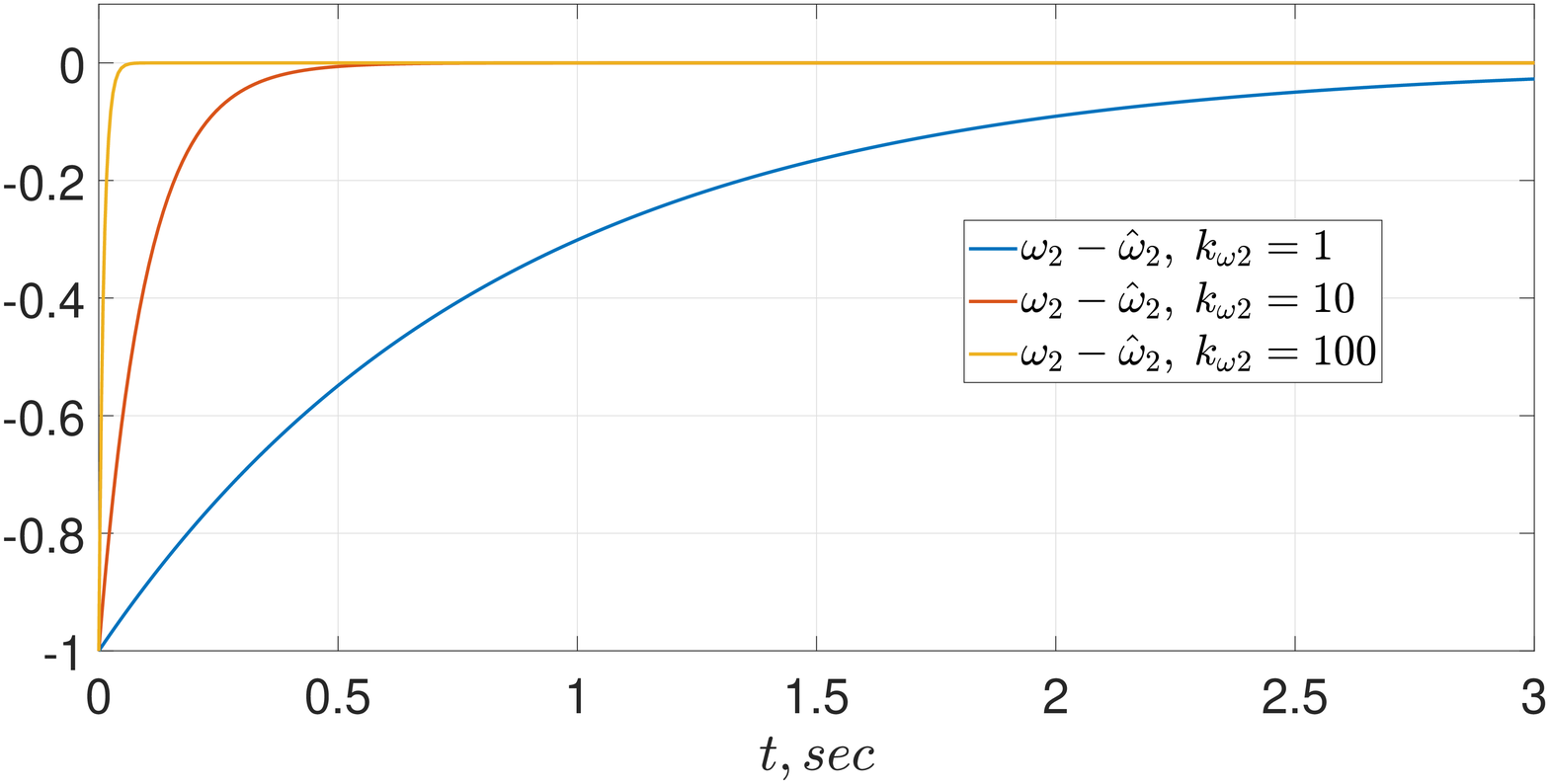}    
		\caption{Transients of the second speed observation error for the observer \eqref{hat_om} for different values of $k_{w1}$ 
		}
		\label{fig:pic6}
	\end{center}
\end{figure}
\subsection{Chemical-biological reactors}
\label{subsec53}
%
We consider reaction systems whose dynamical model is given by \cite[Section 1.5]{BASDOC}
\begali{
	\nonumber
	\dot c & = - u c+Kr(c )+\chi\\
	y  & = \begmat{I_p & \vdots & 0_{p \times d}}c,
	\lab{sysrea}
}
with $c \in \rea_+^n,\chi \in \rea_+^n$, $u \in \rea_+$, $y \in \rea^p$, $r:\rea^n \to \rea_+^q$, $d:=n-p$, $q<n$. It is assumed that $y,u,\chi$ and $K$ are {\em known}. \\
To simplify the notation we partition the vector $c$ as $c = \col(y,x)$, and rewrite \eqref{sysrea} as
\begali{
	\nonumber
	\dot y&=-uy+K_y r(y,x ) + \chi_y\\
	\dot x&=-u x +K_x r(y,x ) + \chi_x.
	\lab{refsys}
}
To simplify the presentation we assume that there are {\em more} measurements than reaction rates, that is, $p \geq q$ and $\rank\{K_y\}=q$.\footnote{See \cite{ORTetaljpc} for a relaxation of this assumption.}
\subsubsection{Solution via GPEBO}
\label{subsubsec531}
The following lemma identifies the mappings $\phi$, $\Lambda$ and $B$ required to satisfy conditions (i) and (ii) of Proposition \ref{pro1}.
\begin{lem}\em
	\lab{lem2}
	Consider the system \eqref{refsys}. The mappings
	\begali{
		\nonumber
		\phi & :=x - K_x K_y^\dagger y \\
		\nonumber
		\Lambda & :=-u\\
		B &:=- K_x K_y^\dagger \chi_y {+\chi_x}
		\lab{philamb}
	}
	where
	$$
	K^\dagger_y:=(K^\top_y K_y)^{-1}K^\top_y,
	$$  
	satisfy the PDE \eqref{pde}. More precisely,
	\begequ
	\lab{dotphi1}
	\dot \phi=\Lambda \phi +B.
	\endequ
\end{lem}
\begin{pf}
	From  \eqref{refsys} and \eqref{philamb} we get
	\begalis{
		\dot \phi&= -u x +K_x r(y,x ) + \chi_x \\
		&- K_x K_y^\dagger[-uy+K_y r(y,x ) + \chi_y]\\
		&=-u \phi + \chi_x- K_x K_y^\dagger \chi_y,
	}
	completing the proof.
\end{pf}
Now, note that from \eqref{dotxi}, \eqref{dotphi} and \eqref{dotphi1} we can, invoking the arguments used in the proof of Proposition \ref{pro1}, establish the relation
\begequ
\lab{phixiphi}
\phi = \xi + \Phi_\Lambda \theta,
\endequ
for some $\theta \in \rea^d$. To obtain a {\em bona fide} regressor equation, that is a linear relation between measurable signals and $\theta$ we would assume condition (iii) of Proposition \ref{pro1}. That is, assume the existence of measurable mappings $C$ and $L$ such that \eqref{algcon} holds, that is $L\phi=C$. Unfortunately, in this example it is not possible to satisfy this condition. However, we can still obtain the required linear regression, needed for the parameter estimation using DREM, as shown in the lemma below.
\begin{lem}\em
	\lab{lem3}
	Assume that the rate vector $r(y,x )$ depends {\em linearly} on the unmeasurable components of the state $x$, that is, it is of the form
	\begequ
	\lab{linrea}
	r(y,x )=R(y )x
	\endequ
	where $R: \rea^p \to \rea^{q \times d}$ is a known matrix.\footnote{See \cite{ORTetaljpc} for the case of nonlinear dependence on $x$.} There exists {\em measurable} signals $\caly \in \rea^d$ and $\Delta \in \rea$ such that
	\begequ
	\lab{scareg}
	\caly = \Delta \theta.
	\endequ
\end{lem}
\begin{pf}
	Defining the partial coordinate $y^\dagger = K_y^\dagger y$, we see from \eqref{refsys} that its dynamics takes the form
	\begali{
		\nonumber
		\dot y^\dagger &= -u y^\dagger + R(y)x + K_y^\dagger\chi_y\\
		\nonumber
		&= -u y^\dagger+ K_y^\dagger\chi_y + R(y)( \phi + K_x y^\dagger)\\
		\nonumber
		&= -u y^\dagger+ K_y^\dagger\chi_y + R(y)( \xi+ \Phi_\Lambda \theta + K_x y^\dagger )\\
		&=  \Psi \theta + \chi_l
		\lab{lre}
	}
	where we used \eqref{philamb} to get the second identity, \eqref{phixiphi} in the third identity and we defined the measurable signals
	\begalis{
		\chi_l &:=  -u y^\dagger+ K_y^\dagger\chi_y + R(y)( \xi + K_x y^\dagger) \\
		\Psi & := R(y)\Phi_\Lambda.
	}
	Applying the filter ${\lambda  \over \bfp+\lambda}$---with $\lambda>0$ a {\em free} tuning parameter---to  \eqref{lre}, and regrouping terms, we obtain the linear regression equation\footnote{As usual in adaptive control, we neglect an additive exponentially decaying term in \eqref{vecreg} that is due to the filters initial conditions.}  
	\begequ
	\lab{vecreg}
	Y=\Psi_f \theta.
	\endequ 
	where we defined the signals 
	\begali{
		\nonumber
		\dot \Psi_f&=-\lambda  \Psi_f +\lambda \Psi\\
		\lab{ydotpsi}
		Y&={\lambda \bfp \over \bfp+\lambda}[y^\dagger] -{\lambda  \over \bfp+\lambda}[\chi_l].
	}
	Multiplying \eqref{vecreg} by $\adj\{ \Psi_f^\top \Psi_f\} \Psi_f^\top$ we obtain the identity \eqref{scareg}, where we defined
	\begali{
		\nonumber
		\caly&:=\adj\{ \Psi_f^\top \Psi_f\} \Psi_f^\top Y \\
		\lab{calydel}
		\Delta&:=\det\{ \Psi_f^\top \Psi_f\},
	}
	This completes the proof.
\end{pf}
\subsubsection{Simulations}
\label{subsubsec532}
To illustrate the performance of the DREM observer proposed in the previous section we consider the model of the anaerobic digestion reactor reported in \cite{marcos2004adaptive}. The dynamics, given in equations (55)-(59) of  \cite{marcos2004adaptive}, {may be} written in the form \eqref{refsys}, \eqref{linrea} with the choices $n=4,q=2,p=2$ 
\begalis{
	K_y &=\begmat{-k_3 & 0 \\ k_4 & -k_1},\;K_x=I_2\\
	R(y)&=\begmat{\mu_1(y_1) & 0 \\ 0 & \mu_2(y_2)},\;{\chi_y}=\begmat{u s_{1,0}\\ u s_{2,0} },\;{\chi_x}=0,
}
where $y_1$, $x_1$, $y_2$ and $x_2$ represent the organic matter concentration (g/l), the acidogenic bacteria concentration (g/l), the volatile fatty acid concentration (mmol), the methanogenic bacteria concentration (g/l) and $u$ is the dilution rate. The positive constants $s_{1,0}$ and $s_{2,0}$ denote the concentration of the substrate in the feed, and $k_1$, $k_3$ and $k_4$  are yield positive coefficients. \\
The two specific growth rates $\mu_1$ and $\mu_2$ are given by
$$
\begmat{\mu_1(y_1)\\  \mu_2(y_2)}=\begmat{\frac{\mu_{m,1} y_1}{K_{S,1}+y_1}\\ \frac{\mu_{m,2} y_2}{K_{S,2}+y_2+K_I y_2^2}}.
$$
where $\mu_{m,1}$, $\mu_{m,2}$, $K_{S,1}$, $K_{S,2}$ and $K_I$ are yield positive coefficients.\\
Notice that $K_y$ is square and full rank, consequently
$$
y^\dagger=K_y^{-1}y=-\begmat{{y_1 \over k_3} \\ {y_2 \over k_1}+ {k_4 y_1 \over k_1 k_3}}.
$$
To design the observer we first identify the signals \eqref{philamb} of Lemma \ref{lem2} as
\begalis{
	\Lambda & =-u\\
	B &=- K_y^{-1} \chi_y=-u\begmat{{ -s_{1,0} \over k_3} \\ {-s_{2,0} \over k_1}- {k_4  s_{1,0} \over k_1 k_3}}.\\
}
Consequently, \eqref{dotxi} and \eqref{dotphi} become
\begalis{
	\dot{ \xi } & = -u  \xi +u\begmat{{ s_{1,0} \over k_3} \\ {s_{2,0} \over k_1}+ {k_4  s_{1,0} \over k_1 k_3}}\\
	\dot \Phi_\Lambda &= -u \Phi_\Lambda,\;\Phi_\Lambda(0)=I_n.
}
Then, we follow the proof of Lemma \ref{lem3} to construct the signals
\begalis{
	\chi_l& =  u \begmat{{y_1 \over k_3} \\ {y_2 \over k_1}+ {k_4 y_1 \over k_1 k_3}}- u\begmat{{s_{1,0} \over k_3} \\ {s_{2,0} \over k_1}+ {k_4  s_{1,0} \over k_1 k_3}}\\ &+\begmat{\mu_1(y_1)[\xi_1 - {y_1 \over k_3}] \\ \mu_2(y_2)[ \xi_2-{y_2 \over k_1}- {k_4 y_1 \over k_1 k_3}]} \\
	\Psi & =\begmat{\mu_1(y_1) & 0 \\ 0 & \mu_2(y_2)} \Phi, 
}
that, together with \eqref{ydotpsi} and \eqref{calydel}, define $\caly$ and $\Delta$ of \eqref{scareg}. The design is completed with the parameter estimator \eqref{gpebo}.\\
For the simulations we used the parameters of \cite{marcos2004adaptive}, that is, $k_1 = 268$  mmol/g, $k_3 = 42.14$, $k_4 = 116.5$ mmol/g, $\alpha = 1$,\footnote{In  \cite{marcos2004adaptive} there is a constant $\alpha=0.5$ entering into the dyamics of $x$ as $\dot x=-\alpha u x +K_x r(y,x ) + \chi_x$. To avoid cluttering the notation, and without loss of generality, we assume this constant is equal to one.} $\mu_{m,1} = 1.2$ $d^{-1}$, $K_{S,1} = 8.85$  g/l, $\mu_{m,2} = 0.74$ $d^{-1}$, $K_{S,2} = 23.2$ mmol, $K_I = 0.0039$ mmol$^{-1}$, $S_{1,0}=1$, $S_{2,0}=1$ and $u=0.1$. \\
The initial conditions for the anaerobic digester were set to $x_1(0)=0.1$ g/l, $y_1(0)=0.05$ g/l, $x_2(0)=0.5$ g and $y_2(0)=4$  mmol/l. We used $\lambda=100$ in the filters of  \eqref{ydotpsi}. Fig.~\ref{fig3} and Fig.~\ref{fig4} show the transient behavior of the state estimation errors for different values of the adaptation gain, with $\gamma=0$ corresponding to the OLO. Notice that, although the convergence rate is increased with larger $\gamma$, an undesirable peak appears {at the beginning of the  error transient}. 
\begin{figure}[h]
	\centering
	\includegraphics[width=1 \linewidth]{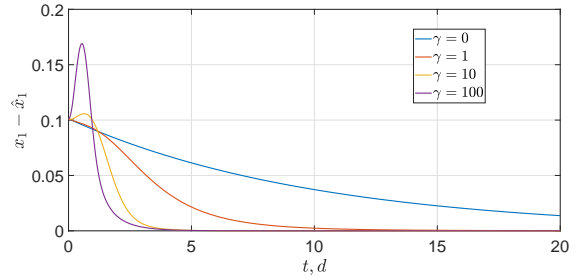}
	\caption{Transients of the error $x_1-\hat{x}_1$ }
	\label{fig3}
\end{figure} 
\begin{figure}[h]
	\centering
\includegraphics[width=1 \linewidth]{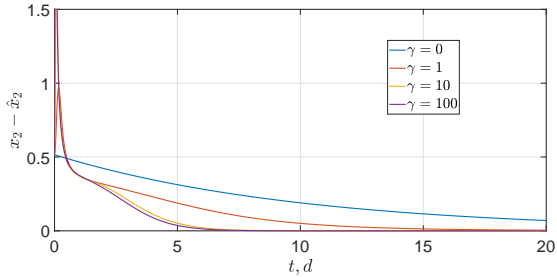}
	\caption{Transients of the error $x_2-\hat{x}_2$ }
	\label{fig4}
\end{figure} 
%
\section{Concluding Remarks}
\label{sec6}
%
An extension to the PEBO technique reported in \cite{ORTetalscl} has been proposed in the paper. It allows us to simplify the task of solving the key PDE and avoid a, sometimes problematic, open-loop integration required in PEBO. Also, we have identified a condition---verification of the algebraic equation \eqref{algcon}---that trivializes the task of estimating the unknown parameters. In the original version of PEBO this was left as an open problem to be solved. It is shown that this condition is satisfied for the practically important problem of power systems.\\
It has been shown that combining PEBO with DREM it is possible, on one hand, to relax the excitation conditions to ensure parameter convergence. On the other hand, it allows us to design an observer with FCT under weak excitation assumptions.  \\
As an additional example we show the application of PEBO+DREM to reaction systems. Notice that the use of DREM is necessary to solve the parameter estimation problem in this example. Although there are many ways to design an estimator from the linear regression \eqref{vecreg}, there exists a fundamental obstacle to ensure its convergence. Indeed, from the definition of $\Phi_\Lambda$, that is $\dot \Phi_\Lambda = -u \Phi_\Lambda$ with $u(t)>0$, we have that $\Phi_\Lambda(t) \to 0$, hence  $\Psi(t) \to 0$---loosing identifiability of the parameter $\theta$.  In particular the matrix $\Psi$ {\em cannot satisfy} the well-known persistency of excitation condition 
$$
\int_t^{t + \kappa} \Psi^\top(s)\Psi(s)ds \geq \kappa I_d,
$$ 
which is the necessary and sufficient condition for exponential convergence of the classical gradient and least-squares estimators  \cite[Theorem 2.5.1]{SASBOD}.

%
\section*{Acknowledgements}
%
This paper is partially supported by the Ministry of Science and Higher Education of Russian Federation, passport of goszadanie no. 2019-0898.

%

\bibliographystyle{plain}        
\bibliography{ortega_automatica20_gpebo_v5_full}           

\appendix
\section{Proof of Proposition \ref{pro1}}
\lab{appa}
%
From \eqref{pde} we have that
$$
\dot \phi = \Lambda \phi+B.
$$
Hence, defining the error signal
\begequ
\lab{e}
e:= \phi-\xi 
\endequ
and taking into account the $\xi$ dynamics of the observer, we obtain an LTV system $\dot e = A(t) e$ where we defined  $A(t):=\Lambda(u(t),y(t))$. Now, from the {\eqref{dotphi}} we see that $\Phi$ is the { {\em fundamental matrix}} of the $e$ system, which is bounded in view of condition (iv). Consequently, there exists a {\em constant} vector $\theta \in \rea^n$ such that
$$
e = \Phi \theta,
$$ 
namely $\theta=e(0)$. We now have the following chain of implications
\begalis{
	e = \Phi \theta & \Leftrightarrow \; \phi=\xi+ \Phi \theta \quad (\Leftarrow \eqref{e}) \\
	& \Rightarrow \; L \phi=L \xi+ L \Phi \theta  \quad (\Leftarrow L \times)\\
	& \Rightarrow \; C - L \xi = L \Phi \theta \quad (\Leftarrow \eqref{algcon})\\
	& \Leftrightarrow \; C - L \xi = \Psi \theta \quad (\Leftarrow \eqref{Psi})\\
	& \Rightarrow \; \Psi^\top  (C - L \xi) = \Psi^\top  \Psi \theta \quad (\Leftarrow  \Psi^\top \times) \\
	& \Rightarrow \; Y = \Omega \theta  \quad \Big(\Leftarrow {\lambda  \over \bfp+\lambda}[\cdot]\; \mbox{and} \; \eqref{doty}, \eqref{dotome}\Big)\\
	& \Rightarrow \; \Delta \theta = \caly, \quad (\Leftarrow \adj\{\Omega\} \times\; \mbox{and} \;  \eqref{caly}, \eqref{del} )
}
where we have used the fact that for any, {\em possibly singular}, $n \times n$ matrix $K$ we have $\adj\{K\}K=\det\{K\}I_n$ in the last line. 

From $\phi=\xi+ \Phi \theta$ and \eqref{phil} it is clear that, if $\theta$ is known, we have that
\begequ
\lab{goodx}
x=\phil(\xi+\Phi \theta,y).
\endequ
Hence, the remaining task is to generate an {\em estimate} for $\theta$, denoted $\hat \theta$, to obtain the observed state via $\hat x=\phil(\xi+\Phi\hat \theta,y)$. This is, precisely, generated with \eqref{gpebo}, whose error equation is of the form  
\begequ
\lab{dottilthe}
\dot {\tilde \theta} =-\gamma \Delta^2 \tilde \theta,
\endequ
where ${\tilde \theta}:= \hat \theta -\theta$. The solution of this equation is given by
\begequ
\lab{tilthe}
{\tilde \theta}(t)=e^{-\gamma \int_{0}^{t} \Delta^2(s)ds}\tilde \theta(0).
\endequ
Given the standing assumption on $\Delta$ we have that $\tilde \theta(t) \to 0$. Hence, invoking \eqref{hatx} and \eqref{goodx} we conclude that  $\tilde x(t) \to 0$, where ${\tilde x}:= \hat x -x$. 

\section{Proof of Proposition \ref{pro2}}
\lab{appb}

First, notice that the definition of $w_c$ ensures that $\hat x$, given in \eqref{ftcgpebo}, is well-defined. Now, from (\ref{tilthe}) and the definition of $w$ we have that 
$$
\tilde \theta =w\tilde \theta(0).
$$ 
Clearly, this is equivalent to
$$
(1 - w)\theta = \hat \theta  - w \hat \theta(0).
$$
On the other hand,  under Assumption \ref{ass1}, we have that $w_c(t)=w(t),\; \forall t \geq t_c$. Consequently, we conclude that 
$$
{1 \over 1 - w_c}[\hat \theta - w_c \hat \theta(0)]=\theta,\; \forall t \geq t_c.
$$
Replacing this identity in \eqref{ftcgpebo} completes the proof.
\section{Parameters of the Power System Example}
\lab{appc}

\begin{table}[h]
	\label{table_1}
	\begin{center}
		\begin{tabular}{|c|c|c|}
			\hline
			Parameter & Initial values & After load change\\
			\hline
			$Y_{12}$ & 0.1032 &  0.1032\\
			\hline
			$Y_{21}$ & 0.1032 & 0.1032\\
			\hline
			$b_1$ &0.0223 & 0.02236\\
			\hline
			$b_2$ &0.0265 & 0.0265\\
			\hline
			$D_1$ &1 & 1\\
			\hline
			$D_2$ &0.2 & 0.2\\
			\hline
			$\nu_1$ &1 & 1\\
			\hline
			$\nu_2$ & 1 & 1\\
			\hline
			$B_{11}$ & -0.4373 & -0.5685\\
			\hline
			$B_{22}$ & -0.4294 & -0.5582\\
			\hline
			$G_{11}$ & 0.0966 & 0.1256 \\
			\hline
			$G_{22}$ &0.0926 &0.1204 \\
			\hline
			$a_1$ & 0.2614 & 0.2898\\
			\hline
			$a_2$ & 0.2532 & 0.2864\\
			\hline
			$P_1$ & 28.22 & 28.22\\
			\hline
			$P_2$  & 28.22 & 28.22\\
			\hline
			$E_{f1}$ & 0.2405 & 0.2405\\
			\hline
			$E_{f2}$ & 0.2405 & 0.2405\\
			\hline
		\end{tabular}
	\end{center}
\end{table}
\end{document}